\documentclass[12pt]{amsart}
\usepackage{amssymb}

\newtheorem{theorem}{Theorem}

\newtheorem{proposition}{Proposition}

\newtheorem{remark}{Remark}
\newtheorem{definition}{Definition}
\newtheorem{lemma}{Lemma}

\newcommand{\Z}{{\mathbb Z}}

\newcommand{\R}{{\mathbb R}}
\newcommand{\C}{{\mathbb C}}

\newcommand{\RP}{{\mathbb RP}}
\newcommand{\CP}{{\mathbb CP}}

\begin{document}

\title[]{Lacunas and local algebraicity of volume functions}

\author{V.A.~Vassiliev}

\date{}

\address{Steklov Mathematical Institute and \newline  
National Research University Higher School of Economics, Moscow}

\email{vva@mi.ras.ru}

\thanks{To appear in Journal of Singularities}

\dedicatory{To the memory of Egbert Brieskorn}



\begin{abstract}
The volume cut off by a hyperplane from a bounded body with smooth boundary in $\R^{2k}$ never is an algebraic function on the space of hyperplanes: for $k$=1 it is the famous lemma XXVIII from Newton's Principia.
Following an analogy of these volume functions with the solutions of hyperbolic PDE's, we study the local version of the same problem: can such a volume function coincide with an algebraic one at least in some domains of the space of hyperplanes, intersecting the body?
We prove some homological and geometric obstructions to this integrability property. Based on these restrictions, we find a family of examples of such ``locally integrable'' bodies in Euclidean spaces.
\end{abstract}

\maketitle

\section{Introduction}

According to an Archimedes' theorem, the volume cut by a plane from a ball in $\R^3$ depends algebraically on the coordinates of the plane. The same is true also for all balls and ellipsoids in all odd-dimensional Euclidean spaces, but no additional examples are known by now. 

On contrary, Newton has proved that for no bounded convex domain with smooth boundary in $\R^2$ the areas cut from it by the affine lines depend algebraically on the coordinates of these lines, see \cite{Newton}, \cite{comment}, \cite{Hooke}, \cite{BH}. V.I.~Arnold \cite{Aprob} has conjectured that similar statements hold also in higher dimensions. The even-dimensional part of this problem was completed in \cite{Vassiliev-2014}: there is no bounded domain (convex or not) with smooth boundary in $\R^{2k}$, for which the volume cut off by a hyperplane is algebraic. The odd-dimensional part of the Arnold's conjecture (stating that the ellipsoids in $\R^{2k+1}$ are unique  bodies with this property) has only partial solutions: several geometric obstructions to the algebraicity of volumes are presented in \cite{APLT}, however it is not clear whether they are sufficient for the proof of the general problem. \smallskip

We study a local version of the same problem: given a body $W \subset \R^N$, can the corresponding volume function coincide with an algebraic one at least in some open subset of the space of all affine subspaces in $\R^N$ intersecting $W$?
We prove some topological and geometric obstructions to this {\em local integrability} property, and find a series of new bodies satisfying it.

There is a deep analogy between this problem and the {\em lacuna problem} in the theory of hyperbolic PDE's developed in \cite{Petrovskii 45}, \cite{Leray}, \cite{ABG 70}, \cite{ABG 73}; for a list of parallel notions see page 138 in \cite{APLT}. Many of our objects and terminology are borrowed from the theory of lacunas.

\subsection{Notation and definitions}

Denote by ${\mathcal P}$ the space of all affine hyperplanes in $\R^N$. It almost coincides with $\RP^N$: the homogeneous coordinates $(a_1: \dots : a_N :b)$ define the hyperplane with the equation 
\begin{equation}
a_1 x_1 + \cdots + a_N x_N + b =0,
\label{homog}
\end{equation}
and $(0: \dots : 0:1)$  is the only point in $\RP^N$ but not in ${\mathcal P}$.

Let $W \subset \R^N$ be a {\em smooth body}, that is, a bounded (not necessarily connected) domain with smooth boundary. It defines a two-valued function $V_W$ on ${\mathcal P}$: its values $V_W(X)$ on a hyperplane $X$ are equal to the volumes of 
intersections of the body $W$ with two halfspaces in $\R^N$ separated by $X$.

The space ${\mathcal P}$ consists of open domains whose points are the hyperplanes transversal to $\partial W$, and the walls between these domains formed by the hyperplanes tangent to it: these walls form the {\em projective dual} hypersurface of $\partial W$. Such an open domain in ${\mathcal P}$ is called a {\em lacuna} if the restriction of the volume functions to this domain coincides with an algebraic function on ${\mathcal P}$, that is, there exists a non-trivial polynomial $F(a_1, \dots, a_N, b, V)$ vanishing in any point $(a_1, \dots, a_N, b, V)$ such that $V$ equals either of  the two volumes cut off from the body $W$ by the hyperplane with the equation (\ref{homog}) from our domain. 
The body $W$ is called {\em algebraically integrable} if all domains of ${\mathcal P}$ are lacunas.
\smallskip

There is a trivial example of a lacuna: it is the domain consisting of hyperplanes not intersecting the body $W$, so that the corresponding volume function is equal identically to a pair of constants in it, 0 and the volume of entire $W$. 
Given a body, does it define nontrivial lacunas in ${\mathcal P}$ (so that the corresponding volume functions are not constant)?

In the case of convex $W \subset \R^{2k}$ and infinitely differentiable $\partial W$ the answer is negative (there is only one non-trivial domain in ${\mathcal P}$, and it is not a lacuna); for $k=1$ it is the Newton's lemma XXVIII. The main result of \cite{Vassiliev-2014} says that for an arbitrary bounded body with $C^\infty$-boundary in $\R^{2k}$ 
{\em all} regular domains in ${\mathcal P}$ cannot be lacunas  simultaneously. 

\section{Obstructions to the integrability}

In this section we assume that  the boundary $\partial W$ of the body $W \subset \R^N$ 
is a smooth component (or a collection of components) of the zero set of an irreducible polynomial with real coefficients.

For any generic real hyperplane $X$, we define an $(N-2)$-dimen\-si\-onal complex manifold, and some collection of elements of its $(N-2)$-dimensional homology group, one of which is given by the manifold $X \cup \partial W$, and the others are called {\em vanishing cycles}. Our main result (Theorem \ref{obstr} below) says that if the intersection index of the first cycle with either of these vanishing cycles is not equal to 0, then the component of ${\mathcal P}$ containing $X$ is not a lacuna. Let us introduce all these objects.

Let $A$ be the zero set in $\C^N$ of the polynomial distinguishing $\partial W$. This set $A$ can have singular points in the imaginary domain. 
Let us fix a Whitney stratification of the algebraic subvariety $A \cup \CP^{N-1}_\infty \subset \CP^N$, where $\CP^N$
is the standard compactification of $\C^N$, and $\CP^{N-1}_\infty$ is the ``infinitely distant'' hyperplane in  it.
An affine hyperplane $X \subset \C^N$ is called {\em generic} if its closure in $\CP^N$ is transversal to this chosen stratification of $A \cup \CP^{N-1}_\infty.$ The set of generic hyperplanes contains a Zariski open subset in the space ${\mathcal P}_\C$ of all complex hyperplanes in $\C^N$. In particular, the real planes in $\R^N$, whose complexifications  are generic, are dense in ${\mathcal P}$. Using the complexifications of real planes, we will consider ${\mathcal P}$ as a subset of ${\mathcal P}_\C$.

Denote by $\mbox{Reg}$ the space of all generic hyperplanes in $\C^N$ , and denote by $\mbox{Reg}_\R$ the set of hyperplanes with real coefficients that are transversal to $\partial W$; in particular $\mbox{Reg}_\R \supset \mbox{Reg} \cap {\mathcal P}$. 
All elements of the difference $\mbox{Reg}_{\R} \setminus (\mbox{Reg} \cap {\mathcal P})$ correspond to real planes whose complexifications are not transversal to the stratified variety $A \cup \C P^{N-1}_\infty$ at some pairs of its complex conjugate imaginary points. The codimension of this difference in ${\mathcal P}$ is at least 2, in particular it does not separate different connected components of $\mbox{Reg} \cap {\mathcal P}$.

The volume function is analytic inside any component of $\mbox{Reg}_\R$.

Given a complex hyperplane $X$ in $\C^N$, denote by $\breve \C^N, \breve X$ and $\breve A$ the sets $\C^N, X$ and $A$ from which all singular points of the hypersurface $A$ are removed.

 Consider the chain of homomorphisms 
\begin{equation}
H_{N} (\breve \C^N, \breve  X \cup \breve A) \to H_{N-1}(\breve X \cup \breve A) \to H_{N-2} (\breve X \cap \breve A) ,
\label{BMV}
\end{equation}
where the first arrow is the usual boundary operator, and the second one is the Mayer-Vietoris differential. 
(All homology groups here and below are with integer coefficients only). 

By the Thom isotopy lemma (see e.g. \cite{GM}), for all $X \in \mbox{Reg}$ the groups of any of three kinds indicated in (\ref{BMV}) are isomorphic to each other; moreover, any path in $\mbox{Reg}$ identifies such groups for the endpoints of the path via the {\em Gauss--Manin connection} (that is, the homological realization of the covering homotopy property over this path). 

Let $X_0 \in \mbox{Reg} \cap {\mathcal P}$ be a generic plane. The group $H_{N} (\breve \C^N, \breve  X_0 \cup \breve A)$
contains two important elements $\Lambda_\pm(X_0)$: the parts of the body $W \subset \R^N$ cut off by the real part of the hyperplane $X_0$ and taken with the canonical (once fixed) orientation of $\R^N$. 
Let $\Delta_\pm(X_0)$ be the images of these elements in the group $H_{N-2} (\breve X_0 \cap \breve A)$
under the composite homomorphism (\ref{BMV}). They are represented by the manifold $X_0 \cap \partial W$ taken with some (opposite) orientations, in particular $\Delta_-(X_0) + \Delta_+(X_0)=0$. 

For any $X \in \mbox{Reg}$ the first and the last groups in (\ref{BMV}) contain also some distinguished sets of elements, called {\em vanishing contours} and {\em vanishing cycles} respectively and defined in the following way.

Let $u$ be a generic point of the hypersurface $\breve A$, that is, a non-singular point of $A$ such that the second fundamental form of $A$ at this point is non-degenerate. Such points are dense in $A$ since $A$ is irreducible and bounds a body in $\R^N$. The set of all hyperplanes tangent to $A$ at points close to $u$ is then a smooth hypersurface in ${\mathcal P}_\C$.  

Let $B$ be a small ball in $\C^N$ centered at our generic point $u \in \breve A$, and $X(u) \subset \C^N$ be the tangent hyperplane of $A$ at $u$.
For any hyperplane $X'(u)$ sufficiently close to $X(u)$ but lying in $\mbox{Reg}$, consider the sequence 
\begin{equation}H_N(B, X'(u) \cup A) \to H_{N-1}((X'(u) \cup A) \cap B ) \to H_{N-2}(X'(u) \cap A \cap B),
\label{MVLoc} 
\end{equation}   
whose maps are defined as in (\ref{BMV}). All three groups in this sequence
are then isomorphic to $\Z$, and both maps in it are the isomorphisms. Denote by $\Lambda(u)$ and $\Delta(u)$ some generators of the first and the last groups in (\ref{MVLoc}) obtained one from another by this composite homomorphism.
Denote by the same letters   $\Lambda(u)$ and $\Delta(u)$
 the images of these elements in the groups $H_N(\breve \C^N, \breve X'(u) \cup \breve A)$ and $H_{N-2}(\breve X'(u) \cap \breve A)$ under the identical embedding. 

An arbitrary path in $\mbox{Reg}$ connecting the points $X'(u)$ and $X_0$ identifies the groups of any of three types (\ref{BMV}) for these hyperplanes, in particular moves the elements $\Lambda(u)$ and $\Delta(u)$ into some two elements of the groups  $H_N(\breve \C^N, \breve X_0 \cup \breve A)$ and $H_{N-2}(\breve X_0 \cap \breve A)$ respectively. All elements of the latter two groups which can be obtained in this way from any choice of a generic point $u$, a path connecting $X$ and $X'(u)$ in $\mbox{Reg}$, and a generator of the group $H_N(B,X'(u) \cup A)$, are called the {\em vanishing contours} and {\em vanishing cycles} respectively.

\begin{theorem}
\label{obstr}
If the domain of {\rm $\mbox{Reg}_\R \subset {\mathcal P}$} containing $X_0$ is a lacuna then the intersection indices $\langle \Delta_+(X_0), \Delta \rangle \equiv - \langle \Delta_-(X_0), \Delta \rangle$ of $(n-2)$-dimensional cycles in the complex $(n-2)$-dimensional manifold $\breve X_0 \cap \breve A$ are equal to 0 for all vanishing cycles $\Delta \in H_{N-2}(\breve X_0 \cap \breve A)$.
\end{theorem}

{\it Proof.} 
The integrals of the holomorphic volume form 
\begin{equation}dx_1 \wedge \dots \wedge x_N
\label{volume}
\end{equation}
along the relative cycles define a linear function on 
the group $H_{N} (\C^N, X \cup A)$, and also on the group $H_N(\breve \C^N, \breve X \cup \breve A)$ for any $X \in {\mathcal P}$. 

Every element $\Lambda$ of the group 
\begin{equation}
H_N(\breve \C^N , \breve X_0 \cup \breve A)
\label{hg}
\end{equation}
 defines a function germ \ $\mbox{Int}(\Lambda)$ \  in a neighborhood of our point $X_0$ in $\mbox{Reg}$: 
its value at any point $X \approx X_0$ is equal to the integral of the form  (\ref{volume}) along the relative cycle $\Lambda(X) \in H_N(\breve \C^N , \breve X \cup \breve A)$, obtained from $\Lambda$ by the Gauss-Manin connection over the paths connecting $X_0$ and $X$ in our neighborhood. By the construction, this function is complex analytic. If $\Lambda$ is one of cycles $\Lambda_+$ or $\Lambda_-$, then the restriction of this function to $\mbox{Reg}_\R$ coincides with the volume function, which also is analytic; therefore the analytic continuations of both functions to entire $\mbox{Reg}$ coincide. If this analytic continuation is infinite-valued then the domain of $\mbox{Reg}_\R$ containing $X_0$ is not a lacuna. 

So we get a linear map \ $\mbox{Int}$ \ from the group (\ref{hg}) to the space of all analytic function germs at the point $X_0 \in {\mathcal P}$. 
Denote by ${\mathfrak H}$ the image of the group (\ref{hg}) under this map (or, equivalently, the group (\ref{hg}) itself factored through the subgroup consisting of all elements defining zero germs). By the construction, ${\mathfrak H}$ is an integer lattice. The group $\pi_1(\mbox{Reg}, X_0)$ acts on the group (\ref{hg}) by monodromy operators, and on ${\mathfrak H}$ by analytic continuations; these actions commute with our epimorphism $\mbox{Int} : H_N(\breve \C^N , \breve X_0 \cup \breve A) \to {\mathfrak H}$.

Now suppose that $\langle \Delta_+(X_0), \Delta \rangle \neq 0$ for some cycle $\Delta$ vanishing along a path connecting the points $X_0$ and $X'(u)$. Consider the loop in $\pi_1(\mbox{Reg},X_0)$ going along this path from $X_0$ to $X'(u)$, rotating around the set of planes tangent to $A$ at points close to $u$, and coming back to $X_0$ along the same path. By the Picard--Lefschetz formula (and the functoriality of the maps (\ref{BMV})) this loop adds to the cycle $\Lambda_+(X_0)$ the class of the contour $\Lambda$ vanishing along our path and taken with a non-zero coefficient $c$ (equal to $\pm \langle \Delta_+(X_0), \Delta \rangle$).

If $N$ is odd then we will pass this loop again and again. In this case the intersection index of $(N-2)$-dimensional cycles in $\breve X \cap \breve A$ is skew-symmetric, therefore any new travel along this loop adds to our integration chain a new copy of the cycle $c \cdot \Lambda$. The function germ defined by any vanishing cycle is not equal to zero, hence we get immediately an infinite number of leaves of the analytic continuation.

\begin{lemma} 
\label{lem1}
Let  $N$ be even, then the orbit of the germ defined by any vanishing contour $\Lambda$
under our $\pi_1(\mbox{Reg}, X_0)$-action in ${\mathfrak H}$  is infinite.
\end{lemma}

{\it Proof} of this lemma is based on considerations of \S 3 in \cite{Vassiliev-2014}. The main tool there is a reflection group associated with any body like $W$. It acts on a lattice ${\mathfrak F}$ generated by finitely many elements corresponding to the vanishing contours, and the orbits of all these generators are not greater than the orbit of an arbitrary germ $\mbox{Int}(\Lambda)$ defined by our vanishing contour under the action of the entire group $\pi_1(\mbox{Reg}, X_0)$. (The action by reflections in ${\mathfrak F}$ is defined by the loops in \ $\mbox{Reg}$, \ all whose points are the planes parallel to $X_0$). Therefore if our $\pi_1(\mbox{Reg}, X_0)$-orbit in ${\mathfrak H}$ of a germ defined by a vanishing contour is finite, then this reflection group also should be finite. However, it was proved in   \cite{Vassiliev-2014} that this reflection group always is infinite. \hfill $\Box$ \medskip

Therefore the orbit of our germ  $\mbox{Int}(c \cdot \Lambda)$ also is infinite. However, this orbit is a subset of the set of differences between the elements of the orbit ot the class $\mbox{Int} (\Lambda_+(X)) \in {\mathfrak H}$. The latter orbit is thus also infinite, that is, the analytic continuation of the volume function has infinitely many leaves at the point $X_0$, and cannot be algebraic. \hfill $\Box$ \medskip   

\begin{theorem}
\label{neigh}
If $N$ is even then two {\em neighboring} domains of the set {\rm $\mbox{Reg}_\R$} of generic hyperplanes in ${\mathcal P}$ $($that is, two domains separated by only one piece of the variety projective dual to $\partial W)$ cannot be lacunas simultaneously.
\end{theorem}

{\it Proof.} 
Let $X_1$, $X_2$ be two points of $\mbox{Reg} \cap {\mathcal P}$ separated by such a piece consisting of hyperplanes tangent to the surface $\partial W$ close to some its \underline{generic}  point $u$; suppose that the planes $X_1$ and $X_2$ are parallel and very close to the plane $X(u)$ tangent to $A$ at this point. Then we have three important elements of the group $H_N(\breve \C^N, \breve X_1 \cup \breve A)$. The first one is our real contour $\Lambda_+(X_1)$ defined by the points of $W$ cut off by the plane $X_1$. The second cycle, $M(\Lambda_+(X_2)),$ is obtained from the similar element $\Lambda_+(X_2)$ of the group  $H_N(\breve \C^N, \breve X_2 \cup \breve A)$ by the Gauss--Manin continuation over a small arc connecting the points $X_2$ and $X_1$ in the space $\mbox{Reg}$ of generic complex hyperplanes. The third element is the vanishing cycle $\Lambda(u)$ generating the group $ H_N(B, X_1 \cup A)$ where $B$ is a small ball centered at the point $u$, see (\ref{MVLoc}). 
By Lemma 3.3 of \S III.3 in \cite{APLT}, these three cycles are related by the equality 
\begin{equation}
\Lambda_+(X_1) - M(\Lambda_+(X_2)) = \pm \Lambda(u), 
\label{induc}
\end{equation}
where the sign $\pm$ depends on the choice of the orientation of the last cycle. By Lemma \ref{lem1}, the orbit of the class $\mbox{Int}(\Lambda(u)) \in {\mathfrak H}$ of the vanishing contour $\Lambda(u)$ under the monodromy action in ${\mathfrak H}$ is infinite in the case of even $N$, therefore the orbits of the classes of elements $\Lambda_+(X_1)$ and $\Lambda_+(X_2)$ cannot be finite simultaneously. \hfill $\Box$ 

\begin{remark} \rm
It follows by induction from the identity (\ref{induc}) that either of the relative homology classes $\Lambda_+(X_0)$ and $\Lambda_-(X_0)$ is equal to the sum of several vanishing contours corresponding to the tangency points of $\partial W$ with the hyperplanes parallel to $X_0$ and lying to the corresponding side from it.
\end{remark}

\section{Local geometry of the boundaries of lacunas and  Davydova condition}

Let $X_1$ and $u$ be the same as in the previous proof.
Let $\Delta_+(X_1)$ and $\Delta(u)$ be two elements of the group $H_{N-2}(\breve X_1 \cap \breve A)$ obtained by the homomorphism (\ref{BMV}) from the elements $\Lambda_+(X_1)$ and $\Lambda(u)$ used in this proof. If their intersection index in $\breve X_1 \cap \breve A$ is not equal to zero, then by Theorem \ref{obstr} the domain of $\mbox{Reg}_\R$ containing $X_1$ is not a lacuna. This property $\langle \Delta_+(X_1), \Delta(u)\rangle \neq 0$ can be checked directly in the terms of the local geometry of $\partial W$ at the point $u$: more precisely, in the terms of its second fundamental form, cf. \cite{Davydova 45}, \cite{ABG 73}.

Let us choose affine coordinates $y_1, \dots, y_N$ in $\R^N$ with the origin at the point $u$ in such a way that $y_1=0$ on the tangent hyperplane $X(u)$, and $y_1 >0$ on the examined hyperplane $X_1$ in our neighborhood $B$ of the point $u$.
The hypersurface $\partial W$ is then defined  by an equation of the form $y_1 = \chi(y_2, \dots, y_N)$ in a vicinity of the point $u$. The function $\chi$ is smooth and has a critical point at the origin: $d \chi(0)=0$.  This critical point  is Morse since $u$ is generic.
        
\begin{proposition}[see e.g. \cite{Leray} or Theorem 3.1 in page 183 of \cite{APLT}]
\label{pro1}
$\langle \Delta_+(X_1), \Delta(u)\rangle = 0$ if and only if the {\em positive} inertia index of the quadratic part of the 
Taylor expansion of the function $\chi$ at the critical point is even. 
\end{proposition}

The trivial example occurs when this inertia index is equal to 0: in this case the cycle $\Delta_+(X_1)$ (consisting of all real points of $X_1 \cap A$) is empty close to $u$ and certainly cannot intersect the vanishing cycle $\Delta(u)$ concentrated in the neighborhood of $u$. 

\begin{remark} \rm
This geometric condition is completely analogous to the {\em Davydova condition} in the theory of hyperbolic PDE's, see \cite{Davydova 45},
although the integration cycles and forms in this theory are different. In both theories, the homology classes of the varieties like $X \cap A$ play the crucial role. However, in our case these cycles are related with the $N$-dimensional integration contours by the maps $($\ref{BMV}$)$, while in the hyperbolic science the integration contours lie in some groups similar to our $H_N(\C^N \setminus (X \cup A))$, which in the case of generic $X$ are related to the group $H_{N-2}(X \cap A)$ by the double Leray tube operation.
\end{remark}

Now let $U$ be a connected component of the space $\mbox{Reg}_\R \subset {\mathcal P}$, and $Y \in \partial U$ a hyperplane tangent to $\partial W$. 

\begin{definition}[cf. \cite{ABG 73}] \rm
The domain $U$ is a {\em local lacuna} at the point $Y$ if the volume function $V_W$ coincides with a pair of regular analytic single-valued functions in the intersection of the domain $U$ with some neighborhood of the point $Y$ in ${\mathcal P}$. 
\end{definition}

\begin{proposition}[cf. \cite{ABG 73}]
1. Let $Y \in {\mathcal P}$ be a hyperplane having a generic tangency with $\partial W$ at some point $u$. A domain of {\rm $\mbox{Reg}_\R$} is a local lacuna 
close to this point $Y$ if and only if the condition $\langle \Delta_+(X_1), \Delta(u)\rangle = 0$ from Proposition \ref{pro1} is satisfied for some $($and then for any$)$ neighboring point $X_1$ of this domain. 

2. If a domain is not a local lacuna at some generic point of its boundary, then it also is not a lacuna.
\end{proposition}

The proof of statement 1 essentially repeats that of a similar statement in \cite{ABG 73}: it follows from the removable singularity theorem. The proof of statement 2 uses additionally Theorem \ref{obstr}. \hfill $\Box$ \medskip

So, in the case of even $N$ exactly one of neighboring domains of $\mbox{Reg}_\R$ at a generic point $Y \in \partial W$ is a local lacuna, and the other is not. 

In the case of odd $N$, either both neighboring domains are local lacunas or both are not.
In particular, if $N$ is odd and the hypersurface $\partial W$ contains the points at which the inertia indices of its second 
fundamental quadratic form are odd, then the body $W$ definitely  is not algebraically integrable. 

The study of geometric restrictions preventing a domain to be a local lacuna at more complicated points of its boundary also is parallel to that for hyperbolic PDE's, see \cite{Gording 77}, \cite{Petr}, \cite{APLT}.

\section{Examples of lacunas} 

Let $m=N-3$, so that $\R^N$ is decomposed into the sum $\R_x^3 \oplus \R_y^m$.

Our easiest example is the tubular $\varepsilon$-neighborhood in $\R^N$ of the unit 2-sphere in $\R_x^3$, that is, the body 
defined by the inequality
\begin{equation}
\left(\sqrt{x_1^2 + x_2^2 + x_3^2} -1\right)^2 + (y_1^2 + \dots + y_{m}^2) \leq \varepsilon^2 \ ,
\label{bounda}
\end{equation}
where $0< \varepsilon <1$. (This equation of its boundary is not polynomial, but is obviously equivalent to a polynomial one of degree 4).

There is a much more general class of examples. Instead of $y_1^2+ \cdots + y_m^2,$ consider
an arbitrary smooth function $\psi: \R_y^m \to \R_+$, invariant under the central symmetries in $\R^m_y$, whose unique critical point is a minimum point at the origin, $\psi(0)=0$, and the entire set $\psi^{-1}([0,\varepsilon^2])$ is contained
in some compact neighborhood of the origin in $\R_y^m$. Define the body $W$ in $\R_x^3 \oplus \R_y^m$ by the condition
\begin{equation}
\left(\sqrt{x_1^2 + x_2^2 + x_3^2}-1\right)^2 + \psi(y_1, \dots, y_m) \leq \varepsilon^2.
\label{bounda3}
\end{equation}

Denote by $C$ the volume of this body (\ref{bounda3}), and by $\Omega$ the $(N-1)$-dimensional Euclidean volume of its section by an arbitrary hyperplane in $\R^{m+3}$ containing the plane $\R_y^m$. 

\begin{theorem}
\label{3thm}
If a hyperplane $X \subset \R^{3+m}$ defined by some equation
$$\alpha_1 x_1 + \alpha_2 x_2 + \alpha_3 x_3 + \sum_{j=1}^m \beta_j y_j = \gamma$$
is sufficiently close to one containing the subspace $\R_y^m$
$($that is, $X$ is nearly orthogonal to $\R_x^3$ and contains a point of $\R_x^3$ sufficiently close to the origin$)$, 
then the volumes of two parts cut by $X$ from the body $($\ref{bounda3}$)$ are equal to
\begin{equation}
\frac{C}{2} \ \pm \ \Omega  \frac{\gamma}{\sqrt{\alpha_1^2 + \alpha_2^2 + \alpha_3^2}}.
\label{3eq}
\end{equation} 
 In particular, the domain in ${\mathcal P}$ containing $X$ is a lacuna.
\end{theorem}

\begin{remark} \rm
The right-hand fraction in $($\ref{3eq}$)$ is the distance from the plane $X \cap \R^3_x$ to the origin. The values $($\ref{3eq}$)$ do not depend on the coefficients $\beta_j$ in the equation of $X$.
\end{remark}

\begin{lemma}
\label{mlem}
In the conditions of Theorem \ref{3thm}, the $(m+2)$-dimensional volume of the intersection $X \cap W$ is equal to $\frac{\Omega}{\cos \alpha(X)}$ where $\alpha(X)$ is the angle between $\R^3_x$ and  the normal vector of $X$.
\end{lemma}

{\it Proof of lemma.} 
For any $y \in \R_y^m$, the preimage of $y$ under the canonical projection $W \to \R_y^m$ is empty if $\psi(y)>\varepsilon^2$; if $\psi(y)<\varepsilon^2$ then it is a spherical layer in $\R_x^3$ between the spheres of radii $R=1+\sqrt{\varepsilon^2-\psi(y)}$ and
 $r=1-\sqrt{\varepsilon^2-\psi(y)}$.
Let $\tilde X$ be the hyperplane in $\R^{3+m}$ containing the subspace $\R_y^m$ and such that the 2-planes $X \cap \R_x^3$ and $\tilde X \cap \R_x^3$ are parallel to one another. The orthogonal projection of $X \cap W$ to $\tilde X$ consists of points $(x,y)$ such that $\psi(y) \leq \varepsilon^2, $ and $x$ belongs to a section of the above-described spherical layer (depending on $y$) by a 2-plane  (depending also on $X$). If $X$ is indeed sufficiently close to a vertical hyperplane containing $\R^m_y$, then for any $y$ with $\psi(y) <\varepsilon^2$ this plane section of the layer is an annulus. The area of this annulus does not depend on the choice of this cutting 2-plane: if the distance of this plane from the origin in $\R^3_x$ is equal to $h<r$, then this area is equal to $\pi \left( \sqrt{R^2-h^2}^2 -\sqrt{r^2-h^2}^2 \right) = \pi (R^2-r^2) =4\pi \sqrt{\varepsilon^2-\psi(y)}$.
So, the $(m+2)$-dimensional volume of the projection of $X \cap W$ to $\tilde X$ is equal to
$$4 \pi \int_{\psi(y) \leq \varepsilon^2} \sqrt{\varepsilon^2-\psi(y)} d y, $$ which does not depend on $X$ and hence is equal to the constant $\Omega$. Further, the orthogonal projection of planes multiplies the volumes by the cosine of the angle between the normals of these planes. \hfill $\Box$     \medskip

{\it Proof of Theorem \ref{3thm}.} Let $X_0$ be the plane parallel to $X$ and passing through the origin in $\R^{3+m}$. 
Both values of the volume function at the point $X_0$ are obviously equal to one another and hence to $\frac{C}{2}$. 
For any  $\lambda \in [0, \mbox{dist}(X_0,X)]$ denote by  $X(\lambda)$ the plane obtained from $X_0$ by the parallel shift towards $X$ by
the distance $\lambda$. The derivatives of the volume functions $V_W(X(\lambda))$ over the parameter $\lambda$ are then equal to $\pm$ the volume from Lemma \ref{mlem}. So, when we come to $X$, these volumes grow/decrease by $$\frac{\Omega}{\cos \alpha(X)} \times \mbox{dist}(X_0,X).$$  

Consider the right triangle in $\R^{3+m}$ whose vertices are the origin and its projections to the planes $X$ and $X \cap \R^3_{x}$. Its angle at the origin is equal to $\alpha(X)$, the leg at this vertex is equal to \ $\mbox{dist}(X_0,X)$, and the hypotenuse is exactly the fraction in (\ref{3eq}). \hfill $\Box$ \medskip

\begin{remark} \rm
We see that a {\em locally} algebraically integrable body in $\R^{N}$ (that is, a body having non-trivial lacunas) does not need to be algebraic itself: in fact, only finite smoothness is demanded on the function $\psi(y_1, \dots, y_m)$ participating in the construction of our examples.
\end{remark}

\end{document}